\documentclass[12pt,leqno]{article}
\usepackage{mathrsfs}

\usepackage{amsfonts}
\usepackage{amsmath, amssymb, amsthm, amsfonts}

\numberwithin{equation}{section}

\def\hf{{\textstyle{1\over2}}}
\def\b{\beta}
\def\d{{\,\rm d}}
\def\e{\varepsilon}

\def\f{\varphi}
\def\G{\Gamma}
\def\k{\kappa}
\def\l{\lambda}

\def\={\;=\;}

\def\zt{\zeta(\hf+it)}

\def\D{\Delta}

\def\hf{{\textstyle{1\over2}}}

\def\f{\varphi}

\def\le{\leqslant}
\def\ge{\geqslant}

\begin{document}
\baselineskip=17pt
\title{ On the Dirichlet divisor problem in short intervals  }
\author{Aleksandar Ivi\'c and Wenguang Zhai
\footnote{Wenguang Zhai  is supported by National Natural Science
Foundation of China (Grant No. 11171344) and   Natural Science
Foundation of Beijing (No. 1112010)} }
\date{}

\maketitle
\begin{abstract}
We present several new results involving $\D(x+U)-\D(x)$, where $U = o(x)$
and
$$
\Delta(x):=\sum_{n\leqslant x}d(n)-x\log x-(2\gamma-1)x
$$
is the error term in the classical Dirichlet divisor problem.
\end{abstract}
\section{\bf Introduction}

Define as usual
$$\Delta(x):=\sum_{n\leqslant x}d(n)-x\log x-(2\gamma-1)x,$$
where $d(n) = \sum_{\delta|n}1$ is the sum of all positive divisors of $n$,
and $\gamma = -\G'(1)=0.5772157\ldots\,$ is Euler's constant.
Dirichlet first proved in the 19th century that
 $\Delta(x)=O(x^{1/2}).$
 The exponent $1/2$ was subsequently improved by many authors.
The latest result reads
\begin{equation}
\Delta(x)\;\ll\; x^{131/416}(\log x)^{26947/8320},\; \frac{131}{416}
= 0.314903\ldots,
\end{equation}
which was obtained by  Huxley \cite{Hu}. For $\Delta(x)$ we have the
following well-known conjecture.

{\bf Conjecture 1.} For any $\e>0$
\begin{equation}
\Delta(x) \=  O_\e(x^{1/4+\varepsilon}).
\end{equation}

Conjecture 1 is  supported by the classical mean-square result
\begin{equation}
\int_1^T\Delta^2(x)\d x \;=\;\frac{(\zeta(3/2))^4}{6\pi^2 \zeta(3)}T^{3/2}
+F(T)
\end{equation}
with $F(T)=O_\e(T^{5/4+\varepsilon})$, proved by Cram\'er
\cite{Cr}. His result incidentally also shows that $\D(x) = o(x^{1/4})$
cannot hold as $x\to\infty$. Here and
later $\e\,(>0)$ denotes constants which may be arbitrarily small,
but are not necessarily the same ones at each occurrence, while
$O_{a,b,\ldots}$ means that the implied $O$-constant depends on
$a,b,\ldots\,$. The estimate $F(T)=O_\e(T^{5/4+\varepsilon})$ was
improved to $O(T\log^5 T)$ by Tong \cite{To}, to $O(T\log^4 T)$ by
Preissmann \cite{P}, and recently to $O(T\log^3 T \log\log T)$ by
Lau and Tsang \cite{LT}.

\medskip
Conjecture 1 is also supported by the  upper bound estimate (see
Ivi\'c \cite{Iv1} and \cite{Iv2})
\begin{equation}
\int_1^T|\Delta(x)|^A\d x \;\ll_\e\; T^{1+A/4+\varepsilon},
\end{equation}
where $0\leqslant A\leqslant 35/4.$ The exponent $35/4$ can be replaced by
$262/27$ if we substitute Huxley's exponent $131/416$ into Ivi\'c's
machinery.   For this kind of estimate, we have the following
conjecture.

{\bf Conjecture 2.} The estimate (1.4) holds for any $A>0.$

{\bf Remark 1.} Obviously if Conjecture 1 is true, then Conjecture 2
is also true. It is easy to  show that if Conjecture 2 is true, then
Conjecture 1 is also true. Hence these two conjectures are
equivalent.

For the asymptotic formulae of higher power moments of $\Delta(x)$
see, for example, the papers of Ivi\'c-Sargos \cite{IS},
Tsang \cite{Ts1} and Zhai  \cite{Z}.

\section{\bf Sign changes of $\Delta(x)$ and a result of Jutila}

Suppose $T$ is a large parameter. Ivi\'c \cite{Iv3} proved that
there exists a positive constant $C>0$ such that $\Delta(x)$ changes
its sign on $[T, T+C\sqrt T].$ More precisely, one can find $x_1,
x_2\in [T, T+C\sqrt T]$ such that $\Delta(x_1)>c_1T^{1/4}$ and
$\Delta(x_2)<-c_2T^{1/4}$ hold respectively. This fact was proved
independently in Heath-Brown and Tsang \cite{HBT}.

Heath-Brown and Tsang \cite{HBT} also proved that the above result
is almost best possible. Actually they proved the following theorem.

\medskip
{\bf Theorem A.} {\it In the interval $[T, 2T]$ there are $\gg
T^{1/2}\log^{5} T$ subintervals of length $\gg T^{1/2}\log^{-5} T$
such that on each subinterval one has $|\Delta(x)|\ge c_3 T^{1/4} $
for some $c_3>0.$}

\medskip
In order to prove Theorem A, Heath-Brown and Tsang used a classical
result of Jutila \cite{J} on the divisor problem in short intervals. This is

\medskip

{\bf Theorem B}. {\it If $\,T\ge 2$ and $1\leqslant
U\ll T^{1/2}\ll H\leqslant T,$ then}
\begin{align} \label{eq18}
&\int_T^{T+H}\left(\Delta(x+U)-\Delta(x)\right)^2\d x\\
& \qquad =\frac{1}{4\pi^2}\sum_{n\leqslant
\frac{T}{2U}}\frac{d^2(n)}{n^{3/2}}\int_T^{T+H}x^{1/2}
   {\Bigl|\exp\left(2\pi i (n/x)^{1/2}U\right)-1\Bigr|}^2\d x \nonumber \\
& \qquad +O_\e(T^{1+\varepsilon}+HU^{1/2}T^{\varepsilon}).\nonumber
\end{align}

\medskip
Note that the terminology ``divisor problem in short intervals'' refers
to the fact that in (2.1) we have $\Delta(x+U)-\Delta(x)$ with
$U = o(T)$ as $T\to\infty$, hence the interval $[x, x+U]$ is ``short''.
In the case when $H=T$,  the first author \cite{Iv5} sharpened (2.1)
to an explicit asymptotic formula.
  The term $T^{1+\varepsilon}$ in  (2.1) can be
replaced by $ T\log^4 T $ if we use the method of Preissmann
\cite{P} in conjunction with the proof of Jutila \cite{J}.

We have the
well-known asymptotic formula (see e.g., [8, Chapter 4])
\begin{equation}
\sum_{n\leqslant x}d^2(n)\=xP(\log x)+O_\e(x^{1/2+\varepsilon}),
\end{equation}
where $P(t)$ is a suitable polynomial of degree three in $t.$ Hence from (2.2) and
(2.1) with $T^{1+\varepsilon}$ replaced by $ T\log^4 T$ one gets, for
$1\leqslant U\leqslant T^{1/2}/2\ll H\leqslant T$, that
\begin{equation}
\int_T^{T+H}\left(\Delta(x+U)-\Delta(x)\right)^2\d x\ll
HU\log^3\frac{\sqrt T}{U}+T\log^4 T.
\end{equation}

With the help of (2.3), Heath-Brown and Tsang   proved the following
Lemma 2.1, which combined with  the results on the moments of
$\Delta(x)$ gives Theorem A.

{\bf Lemma 2.1. } Suppose that  $2\leqslant U\leqslant T^{1/2}.$ Then
\begin{equation}
\int_T^{2T} \max_{0\leqslant u\leqslant U}\Bigl|\Delta(x+u)-\Delta(x)\Bigr|^2\d x\ll
TU\log^5 T.\end{equation}

{\bf Remark 2.} Heath-Brown and Tsang did not prove Lemma 2.1 for
$\Delta(x)$ directly. Actually they proved Lemma 2.1 with
$\Delta(x)$ replaced by
$$
E(T) := \int_0^T|\zt|^2\d t - T\bigl(\log(T/(2\pi)) + 2\gamma-1)\bigr),
$$
which represents the error term in the
asymptotic formula of $\zeta(s)$ on the ``critical line'' $\Re s=1/2.$
However the proof for $\Delta(x)$ is very similar,  even a little simpler.

\medskip
The formula (2.1) led Jutila \cite{J} to propose

{\bf Conjecture 3}. For any $0 < \e < 1/4$ and $x^\e \leqslant U \leqslant x^{1/2-\e}$
we have
$$
\D(x+U) - \D(x) \;\ll_\e\; x^\e\sqrt{U}.
\leqno(2.5)
$$
This conjecture is much stronger than the unconditional estimate
$$
\D(x+U) - \D(x) \;\ll_\e\; x^\e{U}\qquad(1 \ll U \leqslant x),\leqno(2.6)
$$
which easily follows from the definition of $\D(x)$ and the elementary
bound $d(n) \ll_\e n^\e$. It is curious that (2.6) has not been proved
yet by the use of  Vorono{\"\i}'s explicit formula for $\D(x)$
(see e.g., [8, Chapter 3]). On the other hand,
from (1.1) one obtains by trivial estimation
$$
\D(x+U) - \D(x) \;\ll_\e\; x^{131/416+\e}\qquad(1\ll U \ll x).\leqno(2.7)
$$

\section{ \bf Bounds for $\Delta(x+U)-\D(x)$}

In this section we shall present new results on the estimation of the difference
$\Delta(x+U)-\D(x)$, both pointwise and in the statistical sense, by
giving an upper bound for the occurrence of large values. Our pointwise bounds
are obtained without the use of the sophisticated exponential sum techniques
which lead to (1.1).

\medskip
{\bf Theorem 1}. {\it We have
$$
\Delta(x+U)-\D(x) \ll_\e x^{1/4+\e}U^{1/4}\qquad(1\ll U \ll x^{3/5}),\leqno(3.1)
$$
and
$$
\Delta(x+U)-\D(x) \ll_\e x^{2/9+\e}U^{1/3}\qquad(1\ll U \ll x^{2/3}).\leqno(3.2)
$$
Moreover, suppose that
$$
|\D(x_r+U)-\D(x_r)| \ge V \gg U^{1/2} \;(\gg 1) \qquad(r = 1,\ldots\,, R-1),\leqno(3.3)
$$
where $X/2 \le x_1 < \ldots < x_R\le X, |x_r-x_s|\ge V$ for $r\ne s$.
If $(\k,\lambda)$ is an exponent pair for which $\k\ne0$, then
for $X^{\l-\k} \le V^{3+2\l-2\k}U^{-2}$ we have
$$
R \ll_\e X^\e\left(XV^{-5}U^2 + X^{(\k+\lambda)/\k}U^{(2\k+2)/\k}V^{-(3+4\k+2\lambda)/\k}\right).
\leqno(3.4)
$$
}

\medskip
{\bf Corollary 3.1}. If we take in (3.4) the exponent pairs
$$
(\k,\lambda) = (1/2,1/2), (2/7,4/7), (1/6, 4/6),
$$
we obtain
$$
R \ll_\e X^\e\left(XV^{-5}U^2 + X^2U^6V^{-12}\right)\qquad(V\ge U^{2/3}),
$$
$$
R \ll_\e X^\e\left(XV^{-5}U^2 + X^3U^9V^{-37/2}\right)\qquad(V \ge U^{14/25}X^{225}),
$$
$$
R \ll_\e X^\e\left(XV^{-5}U^2 + X^5U^{14}V^{-30}\right)\qquad(V \ge U^{1/2}X^{18}).
$$
respectively.

\medskip
For the definition and properties of (one-dimensional) exponent pairs, see
Ivi\'c \cite{Iv2} or Graham-Kolesnik \cite{GK}.
If in the above estimates one could discard the second
term and retain only the term $XV^{-5}U^2$, this would imply
$$
\int_1^X(\D(x+U)-\D(x))^4\d x \ll_\e X^{1+\e}U^2
$$
in a suitable range for $U$, which is a conjecture of M. Jutila \cite{J}.

\begin{proof}
For the proof of Theorem 1
we need the following well-known truncated form of
the Vorono{\"\i} formula for $\D(x)$ (see e.g., [8, Chapter 3]).

 {\bf Lemma 3.1.} For $1\ll N\ll x$ we have
$$
\Delta(x)=\frac{x^{1/4}}{\sqrt 2\pi}\sum_{n\le
N}d(n)n^{-3/4}\cos(4\pi\sqrt{nx}-\pi/4)+O_\e(x^{1/2+\varepsilon}N^{-1/2}).
\leqno(3.5)
$$

\medskip
Thus setting $f(x) := x^{1/4}\cos(4\pi\sqrt{nx}-\pi/4)$ we have from (3.5)
\begin{eqnarray}
\D(x+U)-\D(x)&&=\frac{1}{\pi\sqrt{2}}\sum_{n\le N}d(n)n^{-3/4}\int_x^{x+U}f'(v)\d v
+O_\e(x^{1/2+\varepsilon}N^{-1/2})\nonumber\\
&&= \frac{1}{4\pi\sqrt{2}}\int_x^{x+U}v^{-3/4}\sum_{n\le N}d(n)n^{-3/4}
\cos(4\pi\sqrt{nv}-\pi/4)\d v\nonumber\\
&&- \sqrt{2}\int_x^{x+U}v^{-1/4}\sum_{n\le
N}d(n)n^{-1/4}\sin(4\pi\sqrt{nv}-\pi/4)\d v \nonumber\\&&
+O_\e(x^{1/2+\varepsilon}N^{-1/2}) = \frac{1}{4\pi\sqrt{2}}I_1 -
\sqrt{2}I_2 + O_\e(x^{1/2+\varepsilon}N^{-1/2}),\nonumber
\end{eqnarray}
say. Also note that $I_1$ and $I_2$ are similar in structure, but $I_1$ is of a lower
order of magnitude, so it suffices to estimate $I_2$. By H\"older's inequality for
integrals we have, for $k\in\mathbb N$,
$$
I_2 \ll \max_K \frac{\log x}{x^{1/4}}U^{1-1/k}{\left(\int_x^{x+U}\Bigl|\sum_{K<n\le K'\le2K}
\frac{d(n)}{n^{1/4}}\exp(4\pi i \sqrt{nx})\Bigr|^{k}\d v\right)}^{1/k},\leqno(3.6)
$$
where the maximum is taken over $O(\log x)$ values of $K \ll N$. We
shall use (3.6) with $k=2$ and $k=4$ to obtain (3.1) and (3.2),
respectively.

When $k=2$ the integral in (3.6) equals
\begin{eqnarray}
{}&&\int_x^{x+U}\sum_{K<m,n\le K'}d(m)d(n)(mn)^{-1/4}
\exp\Bigl(4\pi i(\sqrt{m}-\sqrt{n})\sqrt{v}\Bigr)\d v
\nonumber\\
&&
\ll U\sum_{K<n\le 2K}d^2(n)n^{-1/2} + \sum_{K<m\ne n \le2K}
\frac{d(m)d(n)x^{1/2}}{(mn)^{1/4} |\sqrt{m}-\sqrt{n}|}\nonumber\\
&&
\ll_\e UK^{1/2}\log^3K +K^{1+\e}x^{1/2},\nonumber
\end{eqnarray}
where we used trivial estimation for the terms with $m=n$, and otherwise
the standard first derivative test (see e.g., Lemma 2.1 of \cite{Iv2}).
Therefore we obtain
\begin{eqnarray}
{}&&I_2\ll_\e U^{1/2}x^{\e-1/4}(U^{1/2}N^{1/4} + N^{1/2}x^{1/4})\nonumber\\
&& \quad\;\ll_\e x^\e(UN^{1/4}x^{-1/4} + U^{1/2}N^{1/2}).\nonumber
\end{eqnarray}
This gives
$$
\D(x+U)-\D(x) \ll_\e x^\e(UN^{1/4}x^{-1/4} + U^{1/2}N^{1/2} + x^{1/2}N^{-1/2}).
$$
If we choose $N = (x/U)^{1/2}$, then it follows that
$$
\D(x+U)-\D(x) \ll_\e x^\e(x^{1/4}U^{1/4} + U^{7/8}x^{-1/8}) \ll_\e x^{1/4+\e}U^{1/4}
$$
for $1 \ll U \ll x^{3/5}$, as asserted by (3.1). Note that $x^{1/4}U^{1/4}\le U$ for
$U \ge x^{1/3}$, hence for $x^{1/3} \le U \le x^{3/5}$ we obtain an improvement over
(2.6), without the use of exponential sum techniques.

When $k=4$ we use the technique of the proof of (1.12) of
Ivi\'c-Zhai \cite{IZ}, based on an arithmetic result of
Robert-Sargos \cite{RS} involving four square roots, so we omit the
details. The integral in (3.6) is
$$
\ll_\e x^\e U(K^{5/2}x^{1/2}U^{-1} + K).
$$
Hence from (3.6) we infer that
$$
\D(x+U)-\D(x) \ll_\e x^\e(x^{-1/8}N^{5/8}U^{3/4} + x^{-1/4}UN^{1/4} + x^{1/2}N^{-1/2}).
$$
The choice of $N$ this time will be $N = x^{5/9}U^{-2/3}$, valid for $1 \ll U \le x^{5/6}$.
Thus
$$
\D(x+U)-\D(x) \ll_\e x^\e(x^{2/9}U^{1/3} + U^{5/6}x^{-1/9}) \ll_\e x^{2/9+\e}U^{1/3}
$$
for $1 \ll U \le x^{2/3}$, as asserted by (3.2).

\medskip
It remains to prove (3.4). We shall use the method of \cite{Iv1}, also used in Chapter
13 of \cite{Iv2}. From Lemma 3.1 (taking $N = X^{1+\e}V^{-2}$)
and  the condition (3.3) we obtain, for $r = 1,\ldots, R$,
\begin{eqnarray}
{1}&&\ll x^{-1/4}V^{-1}\int_{x_r}^{x_r+U}\Bigl|\sum_{n\le X^{1+\e}V^{-2}}
d(n)n^{-1/4}e(2\sqrt{nv})\Bigr|\d v\nonumber\\
&& \ll x^{-1/4}V^{-1}U\Bigl|\sum_{n\le X^{1+\e}V^{-2}}
d(n)n^{-1/4}e(2\sqrt{nt_r})\Bigr|,\nonumber
\end{eqnarray}
where $e(z) = \exp(2\pi iz)$, and $t_r$ is the point from $[x_r, x_r+U]$ where
the integral above attains its maximum. Hence we may consider the system of points
$$
X/3 \le t_1 < \ldots < t_R \le 4X/3,\quad |t_r-t_s|\gg V\quad(r\ne s)\leqno(3.7)
$$
such that
$$
1 \ll X^{-1/4}V^{-1}U\Bigl|\sum_{n\le X^{1+\e}V^{-2}}
d(n)n^{-1/4}e(2\sqrt{nt_r})\Bigr|\qquad(r = 1,\ldots, R).\leqno(3.8)
$$
Summation of (3.7) over $r$ and an application of the Hal\'asz-Montgomery
inequality (see e.g., the Appendix of \cite{Iv2}) give
\begin{eqnarray}
{R}&&\ll X^{-1/2}V^{-2}U^2\log X\max_{M\le X^{1+\e}V^{-2}}\sum_{r\le R}
\Bigl|\sum_{M<n\le2M}d(n)n^{-1/4}e(2\sqrt{nt_r}\,)\Bigr|^2\nonumber\\
&& \ll_\e  X^{\e-1/2}V^{-2}U^2\max_{M\le X^{1+\e}V^{-2}}
M^{1/2}\max_{s\le R}\Bigl(M + \sum_{r\le R,r\ne s}\frac{M^{1/2}X^{1/2}}{|t_r-t_s|}
\nonumber\\&&
+ \sum_{r\le R}X_0^\k X^{-\k/2}M^{-\k/2+\l}\Bigr)\nonumber\\&&
\ll_\e X^{1+\e}V^{-5}U^{2} + RX_0^\k X^{-1/2-\k/2+\e}V^{-2}U^2X^{\l-\k/2+1/2}V^{-2\l+\k-1}
\nonumber\\&&
\ll_\e X^{1+\e}V^{-5}U^2 + RX_0^\k X^{\l-\k+\e}U^2V^{-2\l+\k-3}.\nonumber
\end{eqnarray}
This in fact holds if $|t_r-t_s|\le X_0$, namely if we estimate the number of points
$R = R_0$, say, in an subinterval of $[X/3, 4X/3]$ of length $\le X_0$ for a given $X_0$
to be determined a little later. Here we used the estimate
$$
\sum_{M<n\le2M}e(2\sqrt{nx}\,) \;\ll\; x^{\k/2}M^{\l-\k/2},
$$
where $(\k, \l)$ is an exponent pair. It follows that
$$
R_0 \;\ll_\e\; X^{1+\e}V^{-5}U^2
$$
provided that
$$
X_0^\k X^{\l-\k+\e}U^2 V^{-2\l+\k-3} \;\ll\;1,
$$
which is satisfied with the choice ($\k >0$)
$$
X_0 \= X^{(\k-\l)/\k+\e}U^{-2/\k}V^{(3-\k+2\l)/\k},\leqno(3.9)
$$
whence
\begin{eqnarray}
{R}&&\ll R_0(1 + X/X_0)
\nonumber\\
&& \ll_\e
X^{1+\e}V^{-5}U^2 + X^{2+\e}V^{-5}U^2X^{(\l-\k)/2}U^{2/\k}V^{(-3+\k-2\l)/\k},\nonumber
\end{eqnarray}
which implies (3.4) of Theorem 1. Since $V \le |t_r - t_s| \le X_0\;(r\ne s)$, we have yet
to check that $V \le X_0$, which is true if $X^{\l-\k} \le V^{3+2\l-2\k}U^{-2}$
and, in view of (3.9), this is the condition given in the formulation of Theorem 1.

\end{proof}
\section{ \bf A new conjecture on $\Delta(x+u)-\D(x)$}

We note that Jutila's result (2.1) holds on the  interval $[T, T+H]$ with
$T^{1/2}\ll H\ll T.$ But Lemma 2.1 is a result on the interval $[T,
2T].$ Comparing (2.3) and (2.4) it is natural to ask if we can find
a short interval type result of (2.4) similar to (2.3). Here we
propose the following Conjecture 3 about this kind of estimate.

{\bf Conjecture 4. } Suppose  $$ \log T\leqslant U\leqslant T^{1/2}/10,\,
 T^{1/2}\ll H\ll T,\, HU\gg T^{1+\varepsilon}.$$  Then the estimate
\begin{equation}
 \int_T^{T+H} \max_{0\leqslant u\leqslant U}\Bigl|\Delta(x+u)-\Delta(x)\Bigr|^2\d x
 \ll HU\log^c T
 \end{equation}
 holds for some absolute constant $c \ge0.$

According to Lemma 2.1, Conjecture 4 is true for $H=T $ with $c=5.$
It is trivially implied by Conjecture 3. Nevertheless, it
 is very strong, since it implies   Conjecture 1. Namely we
 have the following

 {\bf Proposition 4.1.} Conjecture 4 implies Conjecture 1.

\begin{proof}
Suppose $U\ll T^{1/2}\ll H$. Then we have
\begin{eqnarray}
\Delta(T)&&=\frac{1}{U}\int_{T-U}^T
\Delta(T)\d x\\
&&=\frac{1}{U}\int_{T-U}^T \Delta(x)\d x+\frac{1}{U}\int_{T-U}^T
(\Delta(T)-\Delta(x))\d x\nonumber\\
&&\ll \frac{1}{U}\left|\int_{T-U}^T
\Delta(x)\d x\right|+\frac{1}{U}\int_{T-U}^T
\left|\Delta(T)-\Delta(x)\right|\d x\nonumber\\
&&\ll \frac{U+T^{3/4}}{U}+\frac{1}{U}\int_{T-U}^T
\left|\Delta(x+T-x)-\Delta(x)\right|\d x\nonumber\\
&&\ll \frac{U+T^{3/4}}{U}+\frac{1}{U}\int_{T-U}^T \max_{0\leqslant u\leqslant
U}\left|\Delta(x+u)-\Delta(x)\right|\d x\nonumber\\
&&\ll \frac{U+T^{3/4}}{U}+\frac{1}{U}\int_{T-H}^T \max_{0\leqslant u\leqslant
U}\left|\Delta(x+u)-\Delta(x)\right|\d x\nonumber,
\end{eqnarray}
where we used the well-known formula of Vorono{\"\i} \cite{Vo}
\begin{equation}
 \int_0^T\Delta(x)\d x\ = \textstyle\frac{1}{4}T + O(T^{3/4}).
\end{equation}

  By (4.2), the Cauchy-Schwarz inequality and Conjecture 4 with
$U=T^{1/2-\varepsilon}$ and $H=T^{1/2+2\varepsilon} $ we obtain from (4.2)
\begin{eqnarray}
\Delta(T)\ll_\e T^{1/4+\varepsilon}+\frac 1U (HU\log^3
T)^{1/2}H^{1/2}\ll_\e T^{1/4+\varepsilon}.
\end{eqnarray}

 \end{proof}

\section{\bf A partial answer to Conjecture 4}

In this section, we shall show that the argument of Heath-Brown and
Tsang \cite{HBT} implies a partial answer to Conjecture 4. This is

\medskip
{\bf Theorem 2.}
{\it Suppose $\log^2 T\ll U\leqslant T^{1/2}/2,
T^{1/2}\ll H\leqslant T,$ then we have
\begin{align}
 &\int_T^{T+H} \max_{0\leqslant u\leqslant U}\Bigl|\Delta(x+u)-\Delta(x)\Bigr|^2\d x \ll
HU{\cal L}^5+T{\cal L}^4\log{\cal L}\\
&   +H^{1/3}T^{2/3}U^{2/3}{\cal L}^{10/3}(\log {\cal L})^{2/3},\nonumber
\end{align}
where ${\cal L} :=\log T.$}

\begin{proof}

Write $U=2^\lambda b $ where  $\lambda\in {\mathbb N}$ and $1<b\leqslant
U/10$ is a parameter to be determined later. Suppose $v\leqslant u\leqslant 2T.$ By
the definition of $\Delta(x),$ we have
\begin{eqnarray}
\Delta(u)-\Delta(v)&&=\sum_{v<n\leqslant u}d(n)-M(u)+M(v)\\
&&\ge -M(u)+M(v)\nonumber\\
&&\ge -(u-v)(\log u+2\gamma) \nonumber\\
&&\ge  -3(u-v){\cal L},\nonumber
\end{eqnarray}
where $M(z)=z\log z+(2\gamma-1)z.$

Suppose $x\asymp T, 0<u\leqslant U.$ Then there is some integer $j$
such that $0\leqslant j\leqslant U/b$
and  $jb<u\leqslant (j+1)b\leqslant U.$ From (5.2) we have
$$
\Delta(x+jb)-\Delta(x)-3b{\cal L} \leqslant \Delta(x+u)-\Delta(x)\leqslant
\Delta(x+(j+1)b)-\Delta(x)+3b{\cal L} ,
$$
which implies that
\begin{eqnarray}
\max_{0\leqslant u\leqslant U}|\Delta(x+u)-\Delta(x)|&&\leqslant \max_{1\leqslant j\leqslant
2^\lambda}|\Delta(x+jb)-\Delta(x)|+ 3b{\cal L}\\
&&=|\Delta(x+j_0b)-\Delta(x)|+ 3b{\cal L} \nonumber
\end{eqnarray}
for some $1\leqslant j_0=j_0(x)\leqslant 2^\lambda,$ say. We write $j_0$ in the binary system
as
$$
j_0 =2^\lambda\sum_{\mu\in S}2^{-\mu} = 2^{\l-\mu_1} + 2^{\l-\mu_2}
+ \ldots + 2^{\l-\mu_\ell}
$$
for a certain set
$$
S =S(x) = \Bigl\{\; \mu_1, \mu_2, \ldots, \mu_\ell\;\Bigr\},\quad
0 \leqslant \mu_1 < \mu_2 < \ldots < \mu_\ell \leqslant \l
$$ of distinct non-negative integers $\mu_j$. We claim that
\begin{equation}
\Delta(x+j_0b)-\Delta(x)=
\sum_{\mu\in S}\Bigl(\Delta(x+( \nu+1)2^{\lambda-\mu}b)-
\Delta(x+\nu 2^{\lambda-\mu}b )\Bigr),
\end{equation}
where
$$
\nu=\nu_{\mu} = \nu_\mu(x) =\sum_{\rho\in S,\rho<\mu}2^{\mu-\rho}<2^\mu.
$$
The definition of $\nu_\mu$ implies that
$$
\nu_{\mu_1} = 0, \nu_{\mu_2} = 2^{\mu_2-\mu-1}, \ldots,
\nu_{\mu_\ell} = 2^{\mu_\ell-\mu_1} + \ldots + 2^{\mu_\ell-\mu_{\ell-1}}.
$$
Then the right-hand side of (5.4) becomes
$$
\D(x + 2^{\l-\mu_1}b) - \D(x) + \D(x + (2^{\mu_2-\mu_1}+1)2^{\l-\mu_2}b)-
\D(x + 2^{\mu_2-\mu_1}2^{\l-\mu_2}b)
$$
$$
+  \D(x + (2^{\mu_3-\mu_2}+2^{\mu_3-\mu_1}+1)2^{\l-\mu_3}b)
-  \D(x + (2^{\mu_3-\mu_2}+2^{\mu_3-\mu_1})2^{\l-\mu_3}b)+\ldots
$$
$$+ \D(x + (2^{\mu_\ell-\mu_1}+\ldots + 2^{\mu_\ell-\mu_{\ell-1}}+1)2^{\l-\mu_\ell}b)
- \D(x + (2^{\mu_\ell-\mu_1}+\ldots + 2^{\mu_\ell-\mu_{\ell-1}})2^{\l-\mu_\ell}b)
$$
$$
= \D(x + (2^{\l-\mu_1}+\ldots + 2^{\l-\mu_\ell})b)
- \D(x)
= \D(x+j_0b)-\D(x),
$$
since all the other terms cancel out. This establishes (5.4).

\medskip
By the Cauchy-Schwarz inequality we then obtain
\begin{eqnarray*}
&&\ \ \ |\Delta(x+j_0b)-\Delta(x)|^2\\
&&\leqslant |S|\sum_{\mu\in S}\left(\Delta(x+(
\nu+1)2^{\lambda-\mu}b)-\Delta(x+\nu 2^{\lambda-\mu}b )\right)^2.
\end{eqnarray*}

Collecting all possible $\mu $'s and $\nu $'s, we get
\begin{eqnarray}
&&\ \ \ |\Delta(x+j_0b)-\Delta(x)|^2\\
&&\leqslant (\lambda+1)\sum_{\mu\leqslant \lambda}\sum_{0\leqslant \nu<2^\mu}
\left(\Delta(x+( \nu+1)2^{\lambda-\mu}b)-\Delta(x+\nu
2^{\lambda-\mu}b )\right)^2.\nonumber
\end{eqnarray}
Note that now the double sum on the right-hand side of (5.5) is
independent of $x.$ From (5.3) and (5.5) we  immediately see that
\begin{eqnarray*}
&&\ \ \ \max_{0\leqslant u\leqslant U}|\Delta(x+u)-\Delta(x)|^2\\
&&\ll \lambda \sum_{\mu\leqslant \lambda}\sum_{0\leqslant \nu<2^\mu}
\left(\Delta(x+( \nu+1)2^{\lambda-\mu}b)-\Delta(x+\nu
2^{\lambda-\mu}b )\right)^2+b^2{\cal L}^2,
\end{eqnarray*}
which implies that
\begin{eqnarray}
&&\ \ \ \ \ \ \ \ \ \ \ \ \ \ \  \int_T^{T+H}\max_{0\leqslant u\le U}
|\Delta(x+u)-\Delta(x)|^2\d x\\
&&\ll \lambda \sum_{\mu\leqslant \lambda}\sum_{0\leqslant
\nu<2^\mu}\int_T^{T+H} \left(\Delta(x+(
\nu+1)2^{\lambda-\mu}b)-\Delta(x+\nu 2^{\lambda-\mu}b
)\right)^2\d x+Hb^2{\cal L}^2\nonumber\\
&&\ll \lambda \sum_{\mu\leqslant \lambda}\sum_{0\leqslant
\nu<2^\mu}\int_{T+\nu 2^{\lambda-\mu}b}^{T+H+\nu 2^{\lambda-\mu}b}
\left(\Delta(x+ 2^{\lambda-\mu}b)-\Delta(x )\right)^2\d x+Hb^2{\cal
L}^2.\nonumber
\end{eqnarray}

We remark that the error  term $T^{1+\varepsilon}$ in (2.1) can be
replaced by $T{\cal L}^3\log {\cal L}$ if we couple the argument of
Lau and Tsang \cite{LT} with Jutila's proof of (2.1). Hence
similarly to (2.3) we obtain,
 for $1\leqslant U_1\leqslant T_1^{1/2}/2\ll H_1\leqslant T_1$, that
\begin{equation}
\int_{T_1}^{T_1+H_1}\left(\Delta(x+U_1)-\Delta(x)\right)^2\d x\ll
H_1U_1\log^3\frac{\sqrt T_1}{U_1}+T_1\log^3 T_1 \log\log T_1.
\end{equation}

From (5.6) and (5.7), with $T_1 = T + \nu 2^{\l-\mu}b = T + O(U)$
(since $U=2^\lambda b$ and $\nu < 2^\mu$), we infer that
\begin{eqnarray}
&&\ \ \ \ \ \int_T^{T+H}\max_{0\leqslant u\leqslant
U}|\Delta(x+u)-\Delta(x)|^2\d x\\
&&\ll \lambda \sum_{\mu\leqslant \lambda}\sum_{0\leqslant
\nu<2^\mu}(H2^{\lambda-\mu}b{\cal L}^3+T{\cal L}^3\log {\cal
L})+Hb^2{\cal
L}^2\nonumber\\
&&\ll \lambda \sum_{\mu\leqslant \lambda}(H2^\lambda b{\cal L}^3+T2^\mu
{\cal L}^3\log {\cal L})+Hb^2{\cal L}^2\nonumber\\
&&\ll \lambda (H2^\lambda b{\cal L}^4+T2^\lambda {\cal L}^3\log
{\cal L})+Hb^2{\cal L}^2\nonumber\\
&&\ll H2^\lambda b{\cal L}^5+T2^\lambda {\cal L}^4\log {\cal L}
+Hb^2{\cal L}^2\nonumber\\
&&\ll HU{\cal L}^5+TUb^{-1} {\cal L}^4\log {\cal L} +Hb^2{\cal
L}^2.\nonumber
\end{eqnarray}

Now Theorem 2 follows from (5.8) by taking
$$
b\= C\min((TUH^{-1}{\cal
L}^2\log {\cal L})^{1/3}, U/10).
$$
Here $C>0$ is a suitable constant such that one has
$$
\l = \frac{\log U/b}{\log2} \in \Bbb N.
$$
 \end{proof}

From Theorem 2 we get the following Corollary 5.1, which is
well-known but is usually proved by the method of exponential sums.
For a proof of Vorono{\"\i}'s original estimate $\D(x) \ll
x^{1/3}\log x$ without the use of exponential sums, see the first
author's paper \cite{Iv3}.

{\bf Corollary 5.1.} We have the estimate
$$\Delta(x)\ll x^{1/3}\log^{5/3} x (\log\log x)^{1/3}.$$

\begin{proof}
The proof is the same as that of  Theorem 2. But this time we
take $U=T^{1/2}/10, H=10 T^{1/2}.$ We omit the details.
\end{proof}

{\bf Remark 3.}  The whole procedure leading to Corollary 5.1 is as follows:
first we prove the Jutila type result (5.7) from a more accurate form
of Lemma 3.1 (see Meurman \cite{M}), then we prove Theorem 2, and finally we prove
Corollary 5.1. The procedure  begins with Vorono{\"\i}'s formula and is
very long, but the result is a only a little stronger than
$x^{1/3+\varepsilon},$  which is obtained directly from Lemma 3.1 by
taking $N=x^{1/3}.$  So it seems the above procedure is not
interesting.

It is not the case. Note that it is well-known that $\Delta(x)$ has
also the representation
\begin{equation}
\Delta(x) = -2\sum_{n\leq \sqrt x}\psi(x/n)+O(1), \end{equation}
 where
$\psi(t)=\{t\}-1/2$ and $\{t\}$ is the fractional part of $t.$
Actually we can prove (5.7) from (5.9) without using Lemma 3.1,
following the approach given  in Tsang-Zhai \cite{TZ}. And then we prove Theorem 2
and the corollary. This means that we can prove Corollary 5.1 directly
without using Vorono{\"\i}'s formula (Lemma 3.1).

\bigskip

From Theorem 2 we  also get immediately the following Corollary 5.2,
which is a partial answer to Conjecture 4.

{\bf Corollary 5.2.} Suppose that
$$
\log T\leqslant U\leqslant T^{1/2}/10,\;
 T^{1/2}\ll H\ll T.
 $$
Then Conjecture 4 holds for $c=5$ if
$$HU^{1/2}\gg T{\cal L}^{-5/2}\log {\cal L}.$$

\section{\bf Sign changes of $\Delta(x)$ over short intervals}

In this section, we shall give a short interval analogue of Theorem A
via Theorem 2. The result is

\medskip
{\bf Theorem 3.}  {\it Suppose $T, U, H$ are large parameters and $C>1$ is
a large constant such that
$$T^{131/416+\varepsilon}\ll U\leqslant C^{-1}T^{1/2}{\cal L}^{-5},\
\ CT^{1/4}U{\cal L}^5\log {\cal L}\leqslant H\le T.$$
Then in the interval $[T, T+H]$ there are $\gg HU^{-1}$ subintervals
of length   $\gg U$ such that on each subinterval one has $
\pm\Delta(x) \ge c_{\pm} T^{1/4} $ for some $c_{\pm}>0.$}

\medskip
{\bf Corollary 6.1.} Suppose $T, H $ are large parameters  and $C>1$
is a large constant such that $  CT^{3/4} \log {\cal L}\le H\le
T.$ Then in the interval $[T, T+H]$ there are $\gg HT^{-1/2}{\cal
L}^5$ subintervals of length   $\gg T^{1/2}{\cal L}^{-5}$ such that
on each subinterval one has $ \pm\Delta(x) \ge c_{\pm} T^{1/4} $
for some $c_{\pm}>0.$

\begin{proof} We consider only the case
of the $``+"$ sign, and follow the method of proof of Tsang and Zhai
\cite{TZ}. Since $U\gg T^{131/416+\varepsilon}, $ the condition
$H\ge CT^{1/4}U{\cal L}^5\log {\cal L}$ implies $H\gg
T^{235/416+\varepsilon}.$ Thus by Theorem 2 of Lau and Tsang \cite{LT}
we have, as $T\to\infty$,
\begin{equation}
\int_T^{T+H}|\Delta(x)|^2\d x=C_2HT^{1/2}(1+o(1))
\end{equation}
and
\begin{equation}
\int_T^{T+H}|\Delta(x)|^3\d x=C_2HT^{3/4}(1+o(1)),
\end{equation}
where $C_2, C_3$ are suitable positive constants. From (6.1), (6.2) and
the Cauchy-Schwarz inequality we have
\begin{eqnarray*}
 HT^{1/2}&&\ll \int_T^{T+H}|\Delta(x)|^2\d x =
 \int_T^{T+H}|\D^{3/2}(x)\D^{1/2}(x)|^2\d x\\
 &&\ll \left(\int_T^{T+H}|\Delta(x)|^3\d x\right)^{1/2}
 \left(\int_T^{T+H}|\Delta(x)|\d x\right)^{1/2}\\
 &&\ll H^{1/2}T^{3/8}  \left(\int_T^{T+H}|\Delta(x)|\d x\right)^{1/2},
\end{eqnarray*}
which implies that
\begin{equation}
\int_T^{T+H}|\Delta(x)|\d x\gg HT^{1/4}.
\end{equation}
From (4.3) we get
\begin{equation} \int_T^{T+H} \Delta(x) \d x\ll
H+T^{3/4}.
\end{equation}

For any $x\asymp T,$ define
\begin{eqnarray*}
\Delta_{+}(x)=\left\{\begin{array}{ll}
\Delta(x) \quad&\mbox{if $ \Delta(x)>0,$}\\
0\quad& \mbox{otherwise.}
\end{array}\right.
\end{eqnarray*}
We can therefore write
$$
\int_T^{T+H}|\Delta(x)|\d x = 2\int_T^{T+H}\Delta_+(x)\d x - \int_T^{T+H}\Delta(x)\d x.
$$
Then from (6.3), (6.4) and the Cauchy-Schwarz inequality we get
\begin{eqnarray*}
HT^{1/4}&&\ll \int_T^{T+H}|\Delta(x)|\d x \\
&&\ll \left(\int_T^{T+H}\d x\right)^{1/2}\left(\int_T^{T+H}|\Delta_+(x)|^2\d x\right)^{1/2}\\
 &&\ll H^{1/2}   \left(\int_T^{T+H}|\Delta_+(x)|^2\d x\right)^{1/2},
\end{eqnarray*}
which implies that
\begin{equation}
\int_T^{T+H}|\Delta_+(x)|^2\d x\gg HT^{1/2}.
\end{equation}

Finally let us define
$$
\omega(x)=|\Delta_+(x)|^2-4\max_{0\le u\le
U}|\Delta(x+u)-\Delta(x)|^2-\delta x^{1/2},
$$
where $\delta>0$ is a sufficiently small positive constant. If
$\omega(x)>0$, then it follows that
$$\Delta(x) \geq \sqrt{\delta} x^{1/4} $$
and
$$ \Delta(x) \geq 2\max_{0\le u\le
U}|\Delta(x+u)-\Delta(x)|. $$

The second inequality implies that for any $0\le u\le U,$
$$\frac{1}{2}\Delta(x)\le \Delta(x+u)\le \frac 32\Delta(x),$$
namely $\Delta(x+u)$ has the same sign as $\Delta(x).$

 Under the conditions of Theorem 3, from (6.5) and Theorem 2 we see that
\begin{eqnarray}
    &&\ \ \ \ \   \int_T^{T+H}\omega(x)\ dx\\&&
   \gg HT^{1/2}-C_2\delta HT^{1/2}\nonumber\\&&\ \ \   - C_1\Bigl(HU{\cal
L}^5+T{\cal L}^4\log{\cal L}
   +H^{1/3}T^{2/3}U^{2/3}{\cal L}^{10/3}(\log {\cal L})^{2/3}\Bigr)  \nonumber\\
   &&\gg HT^{1/2}\nonumber
\end{eqnarray}
for sufficiently small $\delta$ and some absolute constants $C_1$ and $C_2$.

Let
 $
 \mathscr{S}=\{t\in [T,T+H]: \omega(x)>0\}.$
 By (6.6), H\"older's  inequality and (6.2)   we
get
\begin{eqnarray*}
HT^{1/2}&&\ll \int_T^{2T}\omega(x)\d x\le \int_{
\mathscr{S}}\omega(x)dx\le
\int_{ \mathscr{S}}\Delta_{+}^2(x)\d x\\
&&\le |
\mathscr{S}|^{1/3}\left(\int_T^{2T}|\Delta(x)|^3dx\right)^{2/3} \ll
|\mathscr{S} |^{1/3}H^{2/3}T^{1/2},
\end{eqnarray*}
which implies  $ |\mathscr{S} |\gg   H. $ This completes the proof
of Theorem 3.

\end{proof}

\section{\bf On a problem of Tsang}

In 2010, Tsang wrote a well-written survey paper \cite{Ts2} about
$\Delta(x),$ in which he proposed the following

{\bf Problem}.  Do there exist intervals $[T,\, T+ H], H = T^\beta$
with $\beta>1/4$ such that
\begin{equation}
\int_T^{T+H}|\Delta(x)|\d x \;\ll\; HT^{1/4-\delta}
\end{equation}
 for some small positive $\delta>0$?

 \medskip
Note that, by \cite{Iv3}, for suitable $C>0$ the interval
$[T, T+C\sqrt{T}]$ contains a point $x_0$ where $\D(x)$ changes sign,
hence $\D(x_0) \ll_\e x^\e$. But then, since $\D(x+U)-\D(x)\ll_\e x^\e(U+1)$
for any $U>0$, we have
\begin{eqnarray*}
\int_{x_0}^{x_0+H}|\D(x)|\d x&&= \int_{x_0}^{x_0+H}|\D(x)-\D(x_0)+\D(x_0)|\d x \\
&&\ll_\e \int_{x_0}^{x_0+H}(x_0^\e H+x_0^\e)\d x \ll_\e H^2x_0^\e \le Hx_0^{1/4-\delta}
\end{eqnarray*}
for $H = x_0^\b, 0< \b < 1/4$, provided that $\e$ and $\delta$ are chosen
sufficiently small. This shows why $\b > 1/4$ was assumed by Tsang in connection with (7.1).

There is another easy case of Tsang's problem.
Namely Theorem 1 of Lau and Tsang \cite{LT} implies that if $1/2<\beta<1,$
then we have
\begin{equation}
\int_T^{T+T^\beta}|\Delta(x)|\d x=C_\beta HT^{1/4}(1+o(1))\qquad(T\to\infty)
\end{equation}
for some constant $C_\beta\,(>0)$. The formula (7.2) obviously disproves
(7.1) for $1/2<\beta<1.$ However, it remains to prove or disprove
(7.1) for $1/4\le\beta\le 1/2$, and this is a difficult problem.

In this section we shall show if $\Delta(x)$ could have enough sign
changes, then (7.1) in this range would be true.

We start by  taking $N=T $ in Lemma 3.1. We have
\begin{eqnarray}
\Delta(x)=F(x)+O_\e(T^\varepsilon), \end{eqnarray} where
 $$F(x):=\frac{x^{1/4}}{\sqrt 2\pi}\sum_{n\le T}\frac{d(n)}{n^{3/4}}\cos(4\pi
\sqrt{nx}-\pi/4).$$ Obviously  $F\in C^{\infty}[T/2,3T].$

  \bigskip

Suppose $I \subseteq [T, 2T]$ is any subinterval such that
$\Delta(x)$ changes its sign in $I.$ Then we can find an $x \in I$
such that $F(x)=0$ or at least $|F(x)|\ll_\e T^\varepsilon.$
Further suppose that $\{x_r\}_{r=1}^R$ is a sequence of points such that
$T<x_1<x_2<\cdots<x_R<2T$  and
\begin{eqnarray}
&&F(x_j) \;\ll_\e\; T^\varepsilon,\ j=1,2,\cdots, R,\\
&&|x_i-x_j|\ge H_0,\ 1\le i<j\le R,
\end{eqnarray}
where $1\ll H_0\ll T^{1/2}\ll R\ll T.$

Let $2<H\le H_0/2.$ For each $1\le j\le R,$ we have for $x_j\le
x\le x_j+ H$ that
\begin{eqnarray}
F(x)&&=F(x)-F(x_j)+F(x_j)\\
&&=F(x_j+x-x_j)-F(x_j)+F(x_j)\nonumber\\
&&\ll_\e \max_{0\le h\le H}|F(x_j+h)-F(x_j)|+T^\varepsilon\nonumber\\
&&\ll_\e \max_{0\le h\le H}|F(x+h)-F(x)|+T^\varepsilon.\nonumber
\end{eqnarray}

So by (7.3), (7.6) and Lemma 2.1   we have that
\begin{eqnarray}
&&\ \ \ \ \ \ \sum_{j=1}^R\int_{x_j}^{x_j+H}|F(x)|^2\d x\\
&&\ll_\e \sum_{j=1}^R\int_{x_j}^{x_j+H} \left(\max_{0\leqslant h\le
H}|F(x+h)-F(x)|\right)^2\d x+RT^{2\varepsilon}\nonumber\\
  &&\ll_\e \int_T^{2T}
 \left(\max_{0\le h\leqslant
H}|F(x+h)-F(x)|\right)^2\d x+RT^{2\varepsilon} \nonumber\\
&&\ll_\e HT\log^5 T+RT^{2\varepsilon}.\nonumber
\end{eqnarray}

Formula (7.7) implies that there is some $1\le j_0\le R$ such
that
$$
\int_{x_{j_0}}^{x_{j_0}+H}|\Delta(x)|^2\d x \ll_\e HTR^{-1}\log^5
T+T^{2\varepsilon},
$$
which combined with the  Cauchy-Schwarz inequality yields
\begin{equation}
\int_{x_{j_0}}^{x_{j_0}+H}|\Delta(x)|\d x \ll_\e
HT^{1/2}R^{-1/2}\log^{5/2} T+H^{1/2}T^{\varepsilon}.
\end{equation}

So if we {\it can take}
$$R\gg T^{1/2+\delta},$$ then Tsang's problem is completely solved.

{\bf Remark 4.} Since $\Delta(x)$   has sign changes in the interval
$[T, T+C\sqrt T]$ for some absolute constant $C,$ it is seen that we
can take $R\gg T^{1/2}$  in the above argument. However, this is
still far from solving Tsang's problem.

\medskip
There is another approach to Tsang's problem, which will be briefly
presented now.
By Lemma 2 of Heath-Brown and  Tsang \cite{HBT}, there are long
intervals where $\D(x)$ does not change sign. In particular,
there are intervals of length $\gg \sqrt{T}\log^{-5}T$ in
$[T, 2T]$ where $\D(x)$ does not change sign. The same is true of
$$
F(x) = F_N(x) := \sum_{n\le N}d(n)n^{-3/4}\cos(4\pi\sqrt{nx}-\pi/4),
$$
where $$T \le x \le 2T,\; N = T^{1/2+2\e+2\delta}.$$
\medskip
Then there exists an interval $[X_0-H, X_0+2H]$, where
$F(X_0-H)=0,$ for $H$ satisfying $H\ll \sqrt{T}\log^{-5}T$, where $F_N(x)$ does not
change sign. Let $\f(x) \;(\ge0)$ be a smooth function supported in
$[X_0-H, X_0+2H]$ such that $\f(x) =1$ in $[X_0, X_0 + H]$
and $\f^{(r)}(x) \ll_r H^{-r}$. Then we have
$$
\int_{X_0}^{X_0+H}|\D(x)|\d x \le \int_{X_0-H}^{X_0+2H}\f(x)|\D(x)|\d x
= \left|\int_{X_0-H}^{X_0+2H}\f(x)\D(x)\d x\right|
$$ $$
= {1\over\pi\sqrt{2}}\left|\int_{X_0-H}^{X_0+2H}x^{1/4}\f(x)F_N(x)\d x\right|
+ O(HT^{1/4-\delta}).
$$
We integrate sufficiently many times the last integral by parts,
getting each time the same type of exponential integral,
with a new factor of order
$\ll T^{1/2}H^{-1}n^{-1/2}$ in the $n$-th term in $F_N(x)$.
This means that we may truncate $F_N(x)$
at $N = T^{1+\e}H^{-2}$, or in other words replace $F_N(x)$ by $F_M(x),\,
M = T^{1+\e}H^{-2}$. The point is that, besides the fact that there are no
absolute value signs in the integral, the sum $F_M(x)$ is {\bf shorter}
than $F_N(x)$, which is significant. Also one should be able to use the
fact that $F(X_0-H)=0$ to show that, for $n$ not large, the initial terms
in $F_M(x)$ and $F_N(X_0-H)$ are small. Namely one cannot make use directly
of $F_N(x)$, even for small $n$, and show that their contribution
is $O(HT^{1/4-\delta})$. But the initial terms in $F_M(x)-F_N(X_0-H)$
are small if $x$ is close to $X_0$. This ought to be taken into account
to show that Tsang's conjecture holds true. If we can prove that
$$
\int_{X_0-H}^{X_0+2H}x^{1/4}\f(x)F_N(x)\d x \ll HT^{1/4-\delta}
$$
with $N = M = T^{1+\e}H^{-2}$ and $T^{1/4}
\ll H \ll \sqrt{T}\log^{-5}T$, then we are done.

\vfill\break

\vskip2cm\sc\small
{A. Ivi\'c, Katedra Matematike RGF-a Universiteta u Beogradu, Dju\v sina 7,
11000 Beograd, Serbia. {\tt  ivic@rgf.bg.ac.rs}

\medskip
W. Zhai, Department of Mathematics, China University of Mining and Technology,
Beijing  100083, China. {\tt zhaiwg@hotmail.com }}


\begin{thebibliography}{99}


\bibitem{Cr}H. Cram\'er, \"Uber zwei S\"atze von Herrn G. H. Hardy,
Math. Z. {\bf 15}(1922), 201-210.

\bibitem{GK} S.W. Graham and G. Kolesnik, Van der Corput's method
of exponential sums, LMS Lecture Notes series {\bf126}, Cambridge
University Press, Cambridge, 1991.
\bibitem{HB} D.R. Heath-Brown, The mean value theorem for the Riemann
zeta-function, Mathematika {\bf25}(1978), 177-184.

\bibitem{HB1} D. R. Heath-Brown, The distribution and moments of the error term
in the Dirichlet divisor problem, Acta Arith. {\bf60}(1992), 389-415.



\bibitem{HBT}D. R. Heath-Brown and K.Tsang, Sign changes of $E(t),
\Delta(x)$ and $P(x),$ J. of Number Theory {\bf 49}(1994), 73-83.

\bibitem{Hu} M. N. Huxley,  Exponential sums and Lattice points III, Proc. London Math. Soc.
{\bf 87}(2003), 591-609.

\bibitem{Iv1} Large values of the error term in the divisor problem,
    Inventiones Math.  {\bf71}(1983), 513-520.

\bibitem{Iv2} A. Ivi\'c, The Riemann zeta-function. John Wiley and Sons, New York, 1985.

\bibitem{Iv3} A. Ivi\'c, Large values of certain number-theoretic error terms, Acta Arith.
{\bf56}(1990), 135-159.

\bibitem{Iv4} A. Ivi\'c, The circle and divisor problem,
Bulletin CXXIX de l'Acad\'emie Serbe des Sciences et des
Arts - 2004, Classe des Sciences math\'ematiques et naturelles,
Sciences math\'ematiques No. {\bf29}, pp. 79-83.

\bibitem{Iv5} A. Ivi\'c, On the divisor function and the Riemann
zeta-function in short intervals,
The Ramanujan Journal: Volume {\bf19}, Issue 2 (2009), 207-224.

\bibitem{IS} A. Ivi\'c  and P. Sargos, On the higher power moments
of the error term in the divisor problem, Illinois J. of Math.
{\bf81}(2007), 353-377.

\bibitem{IZ} A. Ivi\'c  and W. Zhai,
 Higher moments of the error term in the divisor problem (in Russian),
Matemati\v ceskie Zametki {\bf88}(2010), 374-383 (= Math. Notes {\bf88}(2010), 338-346).

\bibitem{J} M. Jutila, On the divisor problem for short intervals,
Ann. Univ. Turkuensis Ser. AI {\bf 186}(1984), 23-30.



\bibitem{LT} Y.-K. Lau and K.-M. Tsang,
On the mean square formula of the error term in the Dirichlet
divisor problem, Math. Proc. Camb. Phil. Soc., Vol. {\bf
146}(2009), 277-287.

\bibitem{M} T. Meurman, On the  mean square of the error term in a generalization
of Dirichlet's divisor problem. Acta Arith. {\bf 74}(1996),
351-364.

\bibitem{P} E. Preissman, Sur la moyenne de la fonction z\^eta.
Nagasaka, Kenji (ed.), Analytic number theory and related topics.
Proceedings of the symposium, Tokyo, Japan, November 11-13, 1991.
Singapore: World Scientific, 119-125 (1993).

\bibitem{RS}  O. Robert and P. Sargos,  Three-dimensional
exponential sums with monomials, J. reine angew. Math. {\bf591}(2006), 1-20.


\bibitem{To} K. C. Tong, On divisor problem III, Acta math. Sinica {\bf 6}(1956), 515-541.

\bibitem{Ts1} K.-M. Tsang, Higher-power moments of $\Delta(x), E(t)$ and $P(x)$,
Proc. London Math. Soc.(3) {\bf 65}(1992), 65-84.

\bibitem{Ts2}K.-M. Tsang, Recent progress on the Dirichlet divisor problem and
the mean square of the Riemann zeta-function, Science China
Mathematics, Vol. {\bf53}(2010), 2561-2572.

\bibitem{TZ}K.-M. Tsang and W. Zhai, Sign changes of the error term
in Weyl's law for the Heisenberg manifolds, Transactions of AMS,
 {\bf 364}(2012), Number {\bf 5}, 2647-2666.

\bibitem{Vo} G.F. Vorono{\"\i}, Sur une fonction transcendante et ses
applications \`a la sommation de quelques s\'eries, Ann. \'Ecole Normale
{\bf21}(3)(1904), 207-268 and ibid. {\bf21}(3)(1904), 459-534.


\bibitem{Z} W. Zhai, On higher-power moments of $\Delta(x)$ (II), Acta Arith.
{\bf 114}(2004), 35-54.



\end{thebibliography}
\end{document}